\documentclass[11pt]{article}
\usepackage{epsfig}

\def\be{\begin{equation}}
\def\ee{\end{equation}}
\def\bea{\begin{eqnarray}}
\def\eea{\end{eqnarray}}
\def\bes{\begin{eqnarray*}}
\def\ees{\end{eqnarray*}}

\def\nn{\nonumber}
\def\lb{\label}
\def\bs{\setminus}

\def\T{{\cal T}}
\def\H{{\cal H}}
\def\R{{\bf R}}
\def\C{{\bf C}}
\def\Z{{\bf Z}}

\def\N{{\bf N}}
\def\U{{\bf U}}

\def\Q{{\bf Q}}
\def\T{{\bf T}}

\def\Sg{{\Sigma}}

\def\aa{{\alpha}}
\def\bb{{\beta}}
\def\ga{{\gamma}}

\def\th{{\theta}}

\def\om{{\omega}}
\def\Om{{\Omega}}
\def\ep{{\epsilon}}
\def\lm{{\lambda}}
\def\Lm{{\Lambda}}

\def\sg{{\sigma}}
\def\dm{{\diamond}}
\def\Sg{{\Sigma}}
\def\vf{{\varphi}}

\def\<{{\langle}}
\def\>{{\rangle}}
\def\T{{\cal T}}

\def\Nn{{\cal N}}

\def\mul{{\rm mul}}

\def\crit{{\rm crit}}

\def\Sp{{\rm Sp}}

\def\wtd#1{\widetilde{#1}}

\def\hb{\vrule height0.18cm width0.14cm $\,$}

\title{Non-hyperbolic closed characteristics on non-degenerate star-shaped hypersurfaces in $\R^{2n}$}

\author{Huagui Duan$^{1}$,\thanks{Partially supported by NSFC (Nos. 11131004, 11471169),
LPMC of MOE of China and Nankai University. E-mail: duanhg@nankai.edu.cn.}
\quad Hui Liu$^{2}$,\thanks{Partially supported by NSFC (No. 11401555), China Postdoctoral Science
Foundation No.2014T70589, CUSF (No. WK3470000001). E-mail: huiliu@ustc.edu.cn.}
\quad Yiming Long$^{3}$,\thanks{Partially
supported by NSFC (No. 11131004), MCME and LPMC of MOE of China, Nankai University and BCMIIS of Capital Normal
University. E-mail: longym@nankai.edu.cn.}
\quad Wei Wang$^{4}$ \thanks{Partially supported by NSFC (Nos. 11222105, 11431001).
E-mail: alexanderweiwang@gmail.com.}\\\\
$^{1}$ School of Mathematical Sciences and LPMC, Nankai University, Tianjin 300071\\
$^{2}$ School of Mathematical Sciences, University of Science and Technology of China, Hefei 230026\\
$^{3}$ Chern Institute of Mathematics and LPMC, Nankai University, Tianjin 300071\\
$^{4}$ School of Mathematical Sciences and LMAM, Peking University, Beijing 100871\\
The People's Republic of China\\}

\begin{document}
\maketitle

\begin{abstract}
{\it In this paper, we prove that for every index perfect non-degenerate compact star-shaped hypersurface
$\Sigma\subset\R^{2n}$, there exist at least $n$ non-hyperbolic closed characteristics with even Maslov-type
indices on $\Sigma$ when $n$ is even. When $n$ is odd, there exist at least $n$ closed characteristics with
odd Maslov-type indices on $\Sigma$ and at least $(n-1)$ of them are non-hyperbolic. Here we call a compact star-shaped
hypersurface $\Sigma\subset\R^{2n}$ {\rm index perfect} if it carries only finitely many geometrically distinct prime
closed characteristics, and every prime closed characteristic $(\tau,y)$ on $\Sigma$ possesses positive mean
index and whose Maslov-type index $i(y, m)$ of its $m$-th iterate satisfies $i(y, m)\not= -1$ when $n$ is even,
and $i(y, m)\not\in \{-2,-1,0\}$ when $n$ is odd for all $m\in\N$. }
\end{abstract}



\renewcommand{\theequation}{\thesection.\arabic{equation}}
\renewcommand{\thefigure}{\thesection.\arabic{figure}}

\setcounter{figure}{0}
\setcounter{equation}{0}
\section{Introduction and main results}

Let $\Sigma$ be a $C^3$ compact hypersurface in $\R^{2n}$ strictly star-shaped with respect to the origin, i.e.,
the tangent hyperplane at any $x\in\Sigma$ does not intersect the origin. We denote the set of all such
hypersurfaces by $\H_{st}(2n)$, and denote by $\H_{con}(2n)$ the subset of $\H_{st}(2n)$ which consists of all
strictly convex hypersurfaces. We consider closed characteristics $(\tau, y)$ on $\Sigma$, which are solutions
of the following problem
\be
\left\{\matrix{\dot{y}=JN_\Sigma(y), \cr
               y(\tau)=y(0), \cr }\right. \lb{1.1}\ee
where $J=\left(\matrix{0 &-I_n\cr
        I_n  & 0\cr}\right)$, $I_n$ is the identity matrix in $\R^n$, $\tau>0$, $N_\Sigma(y)$ is the outward
normal vector of $\Sigma$ at $y$ normalized by the condition $N_\Sigma(y)\cdot y=1$. Here $a\cdot b$ denotes
the standard inner product of $a, b\in\R^{2n}$. A closed characteristic $(\tau, y)$ is {\it prime}, if $\tau$
is the minimal period of $y$. Two closed characteristics $(\tau, y)$ and $(\sigma, z)$ are {\it geometrically
distinct}, if $y(\R)\not= z(\R)$. We denote by $\T(\Sigma)$ the set of geometrically distinct
closed characteristics $(\tau, y)$ on $\Sigma\in\mathcal{H}_{st}(2n)$. A closed characteristic
$(\tau,y)$ is {\it non-degenerate} if $1$ is a Floquet multiplier of $y$ of precisely algebraic multiplicity
$2$; {\it hyperbolic} if $1$ is a double Floquet multiplier of it and all the other Floquet multipliers
are not on ${\bf U}=\{z\in {\bf C}\mid |z|=1\}$, i.e., the unit circle in the complex plane; {\it elliptic}
if all the Floquet multipliers of $y$ are on ${\bf U}$. We call a $\Sigma\in \mathcal{H}(2n)$ {\it non-degenerate} if
all the closed characteristics on $\Sigma$ together with all of their iterations are non-degenerate.

The study on closed characteristics in the global sense started in 1978, when the existence of at least one
closed characteristic was first established on any $\Sg\in\H_{st}(2n)$ by P. Rabinowitz in \cite{Rab1}
and on any $\Sg\in\H_{con}(2n)$ by A. Weinstein in \cite{Wei1} independently, since then the existence of
multiple closed characteristics on $\Sg\in\H_{con}(2n)$ has been deeply studied by many mathematicians, for
example, studies in \cite{EkL1}, \cite{EkH1}, \cite{HWZ1}, \cite{Szu1}, \cite{LoZ}, \cite{Wan2}, \cite{Wan3},
\cite{WHL} as well as \cite{Lon4} and references therein.

For the star-shaped hypersurfaces, We are only aware of few papers about the multiplicity of closed
characteristics. In \cite{Gir1} of 1984 and \cite{BLMR} of 1985, $\;^{\#}\T(\Sg)\ge n$ for $\Sg\in\H_{st}(2n)$
was proved under some pinching conditions. In \cite{Vit2} of 1989, C. Viterbo proved a generic existence
result for infinitely many closed characteristics on star-shaped hypersurfaces. In \cite{HuL} of 2002, X. Hu
and Y. Long proved that $\;^{\#}\T(\Sg)\ge 2$ for non-degenerate $\Sg\in \H_{st}(2n)$. In \cite{HWZ2} of 2003,
H. Hofer, K. Wysocki, and E. Zehnder proved that $\,^{\#}\T(\Sg)=2$ or $\infty$ holds for every non-degenerate
$\Sg\in\H_{st}(4)$ provided that all stable and unstable manifolds of the hyperbolic closed characteristics
on $\Sg$ intersect transversally. Recently $\;^{\#}\T(\Sg)\ge 2$ was first proved for every $\Sg\in \H_{st}(4)$
by D. Cristofaro-Gardiner and M. Hutchings in \cite{CGH1} without any pinching or non-degeneracy conditions.
Different proofs of this result can also be found in \cite{GHHM}, \cite{LLo1} and \cite{GiG1}.

On the stability problem, we refer the readers to \cite{Eke1}, \cite{DDE1}, \cite{HuO}, \cite{Lon1}
\cite{Lon3}, \cite{LoZ}, \cite{WHL}, \cite{Wan1} for convex hypersurfaces. For the star-shaped case,
there are also few results. In \cite{LiL2} of 1999, C. Liu and Y. Long proved that for $\Sg\in\H_{st}(2n)$
either there are infinitely many closed characteristics, or there exists at least one non-hyperbolic closed
characteristic, provided every closed characteristic possesses its Maslov-type mean index greater than 2
when $n$ is odd, and greater than 1 when $n$ is even. Recently in \cite{LLo2}, H. Liu and Y. Long proved
that $\Sg\in\H_{st}(4)$ and $\,^{\#}\T(\Sg)=2$ imply that both of the closed characteristics must be elliptic
provided that $\Sg$ is symmetric with respect to the origin.

Based on the common index jump theorem and other results proved by Long and Zhu in \cite{LoZ} of 2002,
recently an enhanced common index jump theorem was established by Duan, Long and Wang in \cite{DLW} and
some new results about the multiplicity and stability of closed geodesics on bumpy simply connected Finsler
manifolds were proved as its applications, which improved some earlier results largely. In this paper, we
will use this enhanced common index jump theorem to study the multiplicity and stability of closed
characteristics on non-degenerate star-shaped hypersurfaces and obtain the following results.

\medskip

{\bf Definition 1.1.} {\it We call a compact star-shaped hypersurface $\Sigma\subset\R^{2n}$ {\rm index
perfect} if it carries only finitely many geometrically distinct prime closed characteristics, and every prime closed
characteristic $(\tau,y)$ on $\Sigma$ possesses positive mean index and whose Maslov-type index $i(y, m)$
of its $m$-th iterate satisfies $i(y, m)\not= -1$ when $n$ is even and $i(y, m)\not\in \{-2,-1,0\}$ when
$n$ is odd for all $m\in\N$.}

Note that our index perfect condition is much more weaker than the dynamically convex condition first
introduced in \cite{HWZ1} (cf. also \cite{HWZ2}), which requires that every closed characteristic on a
compact star-shaped hypersurface $\Sigma\subset\R^{2n}$ possesses Conley-Zehnder index at least $n+1$.
Our theorem below shows that ``index perfect" implies that the number of geometrically distinct closed characteristics
can get to the best possible lower bound at least in the non-degenerate case.

\medskip

{\bf Theorem 1.2.} {\it On every index perfect non-degenerate compact star-shaped hypersurface
$\Sigma\subset\R^{2n}$, there exist at least $n$ non-hyperbolic closed characteristics with
even Maslov-type indices on $\Sigma$ when $n$ is even. When $n$ is odd, there exist at least $n$
closed characteristics with odd Maslov-type indices on $\Sigma$ and at least $(n-1)$ of them are
non-hyperbolic.}

\medskip

{\bf Remark 1.3.} Very recently, J. Gutt and J. Kang in \cite{GuK} proved some results on closed
characteristics on a non-degenerate compact star-shaped hypersurface $\Sg$ in $\R^{2n}$. In Theorem 1.1
of \cite{GuK}, they proved the existence of at least $n$ closed characteristics on such a $\Sg$ whose
iterates' Conley-Zehnder indices possess the same parity, provided $\Sg$ is dynamically convex. Note
that their index definition is slightly different from ours. Also recently, M. Abreu and L. Macarini
in \cite{AbM} gave a sharp lower bound for the number of geometrically distinct contractible periodic
orbits of dynamically convex Reeb flows on prequantizations of symplectic manifolds that are not
aspherical, which implies Theorem 1.1 of \cite{GuK} (cf. Corollary 2.9 of \cite{AbM}). Note that by
\cite{LiL1} and Theorem 1.1 of \cite{LiL3} (cf. Theorem 10.1.2 of \cite{Lon4}), and the definition of
the index in Section 2 of \cite{GuK}, their dynamically convexity implies the mean index of such a
closed characteristic is positive. In Theorem 1.3 of \cite{GuK} a similar result was proved when every
closed characteristic on $\Sg$ possesses Conley-Zehnder index at least $n-1$. But in the proof they
used the common index jump theorem of Y. Long and C. Zhu in \cite{LoZ} and thus they have in fact used
the condition that every closed characteristic possesses positive mean index.

In \cite{LLW2}, H. Liu, Y. Long and W. Wang proved that either there are infinitely many closed
characteristics or there exist at least $2[n/2]$ non-hyperbolic closed characteristics on
$\Sigma\in\H_{con}(2n)$ which is symmetric with respect to the origin, due to the example of weakly
non-resonant ellipsoids, this estimate is sharp when $n$ is even.

Note that comparing with the above three papers, besides non-degenerate condition our above Theorem 1.2
require only that all iterates of closed characteristics do not possess Maslov-type index $-1$ when $n$ is
even, and do not possess Maslov-type indices $-2, -1$ or $0$ when $n$ is odd. We note that the index
conditions in our Theorem 1.2 are much more weaker than the convex or dynamically convex case.

In this paper, let $\N$, $\N_0$, $\Z$, $\Q$, $\R$, $\C$ and $\R^+$ denote the sets of natural integers,
non-negative integers, integers, rational numbers, real numbers, complex numbers and positive real
numbers respectively. We define the function $[a]=\max{\{k\in {\bf Z}\mid k\leq a\}}$, $\{a\}=a-[a]$,
and $E(a)=\min{\{k\in{\bf Z}\mid k\geq a\}}$. Denote by $a\cdot b$ and $|a|$ the standard
inner product and norm in $\R^{2n}$. Denote by $\langle\cdot,\cdot\rangle$ and $\|\cdot\|$
the standard $L^2$ inner product and $L^2$ norm. For an $S^1$-space $X$, we denote by
$X_{S^1}$ the homotopy quotient of $X$ by $S^1$, i.e., $X_{S^1}=S^\infty\times_{S^1}X$,
where $S^\infty$ is the unit sphere in an infinite dimensional {\it complex} Hilbert space.
In this paper we use $\Q$ coefficients for all homological and cohomological modules. By $t\to a^+$, we
mean $t>a$ and $t\to a$.

\setcounter{figure}{0}
\setcounter{equation}{0}
\section{Mean index identities for closed characteristics on compact star-shaped hypersurfaces in $\R^{2n}$}

In this section, we briefly review the mean index identities for
closed characteristics on $\Sg\in\H_{st}(2n)$ developed in \cite{LLW14} which will be needed in Section 4. All
the details of proofs can be found in \cite{LLW14}.
Now we fix a $\Sg\in\H_{st}(2n)$ and assume the following condition on $\T(\Sg)$:

(F) {\bf There exist only finitely many geometrically distinct prime closed characteristics\\
$\qquad\qquad \{(\tau_j, y_j)\}_{1\le j\le k}$ on $\Sigma$. }

Let $\hat{\tau}=\inf_{1\leq j\leq k}{\tau_j}$ and $T$ be a fixed positive constant. Then by Section 2 of
\cite{LLW14}, for any $a>\frac{\hat{\tau}}{T}$, we can construct a function $\varphi_a\in C^{\infty}({\bf R}, {\bf R}^+)$
which has $0$ as its unique critical point in $[0, +\infty)$. Moreover, $\frac{\varphi^{\prime}(t)}{t}$ is strictly
decreasing for $t>0$ together with $\varphi(0)=0=\varphi^{\prime}(0)$ and
$\varphi^{\prime\prime}(0)=1=\lim_{t\rightarrow 0^+}\frac{\varphi^{\prime}(t)}{t}$. More precisely, we
define $\varphi_a$ and the Hamiltonian function $\wtd{H}_a(x)=a\vf_a(j(x))$ via Lemma 2.2 and Lemma 2.4
in \cite{LLW14}. The precise dependence of $\varphi_a$ on $a$ is explained in Remark 2.3 of \cite{LLW14}.

For technical reasons we want to further modify the Hamiltonian, we define the new Hamiltonian
function $H_a$ via Proposition 2.5 of \cite{LLW14} and consider the fixed period problem
\be  \dot{x}(t)=JH_a^\prime(x(t)),\quad x(0)=x(T).  \lb{2.1}\ee
Then $H_a\in C^{3}({\bf R}^{2n} \setminus\{0\},{\bf R})\cap C^{1}({\bf R}^{2n},{\bf R})$.
Solutions of (\ref{2.1}) are $x\equiv 0$ and $x=\rho y(\tau t/T)$ with
$\frac{\vf_a^\prime(\rho)}{\rho}=\frac{\tau}{aT}$, where $(\tau, y)$ is a solution of (\ref{1.1}). In particular,
non-zero solutions of (\ref{2.1}) are in one to one correspondence with solutions of (\ref{1.1}) with period
$\tau<aT$.

For any $a>\frac{\hat{\tau}}{T}$, we can choose some large constant $K=K(a)$ such that
\be H_{a,K}(x) = H_a(x)+\frac{1}{2}K|x|^2   \lb{2.2}\ee
is a strictly convex function, that is,
\be (\nabla H_{a, K}(x)-\nabla H_{a, K}(y), x-y) \geq \frac{\ep}{2}|x-y|^2,  \lb{2.3}\ee
for all $x, y\in {\bf R}^{2n}$, and some positive $\ep$. Let $H_{a,K}^*$ be the Fenchel dual of $H_{a,K}$
defined by
$$  H_{a,K}^\ast (y) = \sup\{x\cdot y-H_{a,K}(x)\;|\; x\in \R^{2n}\}.   $$
The dual action functional on $X=W^{1, 2}({\bf R}/{T {\bf Z}}, {\bf R}^{2n})$ is defined by
\be F_{a,K}(x) = \int_0^T{\left[\frac{1}{2}(J\dot{x}-K x,x)+H_{a,K}^*(-J\dot{x}+K x)\right]dt}. \lb{2.4}\ee
Then $F_{a,K}\in C^{1,1}(X, \R)$ and for $KT\not\in 2\pi{\bf Z}$, $F_{a,K}$ satisfies the
Palais-Smale condition and $x$ is a critical point of $F_{a, K}$ if and only if it is a solution of (\ref{2.1}). Moreover,
$F_{a, K}(x_a)<0$ and it is independent of $K$ for every critical point $x_a\neq 0$ of $F_{a, K}$.

When $KT\notin 2\pi{\bf Z}$, the map $x\mapsto -J\dot{x}+Kx$ is a Hilbert space isomorphism between
$X=W^{1, 2}({\bf R}/({T {\bf Z}}); {\bf R}^{2n})$ and $E=L^{2}({\bf R}/(T {\bf Z}),{\bf R}^{2n})$. We denote its inverse
by $M_K$ and the functional
\be \Psi_{a,K}(u)=\int_0^T{\left[-\frac{1}{2}(M_{K}u, u)+H_{a,K}^*(u)\right]dt}, \qquad \forall\,u\in E. \lb{2.5}\ee
Then $x\in X$ is a critical point of $F_{a,K}$ if and only if $u=-J\dot{x}+Kx$ is a critical point of $\Psi_{a, K}$.

Suppose $u$ is a nonzero critical point of $\Psi_{a, K}$.
Then the formal Hessian of $\Psi_{a, K}$ at $u$ is defined by
\be Q_{a,K}(v)=\int_0^T(-M_K v\cdot v+H_{a,K}^{*\prime\prime}(u)v\cdot v)dt,  \lb{2.6}\ee
which defines an orthogonal splitting $E=E_-\oplus E_0\oplus E_+$ of $E$ into negative, zero and positive subspaces.
The index and nullity of $u$ are defined by $i_K(u)=\dim E_-$ and $\nu_K(u)=\dim E_0$ respectively.
Similarly, we define the index and nullity of $x=M_Ku$ for $F_{a, K}$, we denote them by $i_K(x)$ and
$\nu_K(x)$. Then we have
\be  i_K(u)=i_K(x),\quad \nu_K(u)=\nu_K(x),  \lb{2.7}\ee
which follow from the definitions (\ref{2.4}) and (\ref{2.5}). The following important formula was proved in
Lemma 6.4 of \cite{Vit2}:
\be  i_K(x) = 2n([KT/{2\pi}]+1)+i^v(x) \equiv d(K)+i^v(x),   \lb{2.8}\ee
where the Viterbo index $i^v(x)$ does not depend on K, but only on $H_a$.

By the proof of Proposition 2 of \cite{Vit1}, we have that $v\in E$ belongs to the null space of $Q_{a, K}$
if and only if $z=M_K v$ is a solution of the linearized system
\be  \dot{z}(t) = JH_a''(x(t))z(t).  \lb{2.9}\ee
Thus the nullity in (\ref{2.7}) is independent of $K$, which we denote by $\nu^v(x)\equiv \nu_K(u)= \nu_K(x)$.

By Proposition 2.11 of \cite{LLW14}, the index $i^v(x)$ and nullity $\nu^v(x)$ coincide with those defined for
the Hamiltonian $H(x)=j(x)^\alpha$ for all $x\in\R^{2n}$ and some $\aa\in (1,2)$. Especially
$1\le \nu^v(x)\le 2n-1$ always holds.

For every closed characteristic $(\tau, y)$ on $\Sigma$, let $aT>\tau$ and choose $\vf_a$ as above.
Determine $\rho$ uniquely by $\frac{\vf_a'(\rho)}{\rho}=\frac{\tau}{aT}$. Let $x=\rho y(\frac{\tau t}{T})$.
Then we define the index $i(\tau,y)$ and nullity $\nu(\tau,y)$ of $(\tau,y)$ by
$$ i(\tau,y)=i^v(x), \qquad \nu(\tau,y)=\nu^v(x). $$
Then the mean index of $(\tau,y)$ is defined by
\bea \hat i(\tau,y) = \lim_{m\rightarrow\infty}\frac{i(m\tau,y)}{m}.  \nn\eea
Note that by Proposition 2.11 of \cite{LLW14}, the index and nullity are well defined and are independent of the
choice of $a$. For a closed characteristic $(\tau,y)$ on $\Sigma$, we simply denote by $y^m\equiv(m\tau,y)$
the m-th iteration of $y$ for $m\in\N$.

We have a natural $S^1$-action on $X$ or $E$ defined by
$$  \theta\cdot u(t)=u(\theta+t),\quad\forall\, \theta\in S^1, \, t\in\R.  $$
Clearly both of $F_{a, K}$ and $\Psi_{a, K}$ are $S^1$-invariant. For any $\kappa\in\R$, we denote by
\bea
\Lambda_{a, K}^\kappa &=& \{u\in L^{2}({\bf R}/({T {\bf Z}}); {\bf R}^{2n})\;|\;\Psi_{a,K}(u)\le\kappa\}  \nn\\
X_{a, K}^\kappa &=& \{x\in W^{1, 2}({\bf R}/(T {\bf Z}),{\bf R}^{2n})\;|\;F_{a, K}(x)\le\kappa\}.  \nn\eea
For a critical point $u$ of $\Psi_{a, K}$ and the corresponding $x=M_K u$ of $F_{a, K}$, let
\bea
\Lm_{a,K}(u) &=& \Lm_{a,K}^{\Psi_{a, K}(u)}
   = \{w\in L^{2}(\R/(T\Z), \R^{2n}) \;|\; \Psi_{a, K}(w)\le\Psi_{a,K}(u)\},  \nn\\
X_{a,K}(x) &=& X_{a,K}^{F_{a,K}(x)} = \{y\in W^{1, 2}(\R/(T\Z), \R^{2n}) \;|\; F_{a,K}(y)\le F_{a,K}(x)\}. \nn\eea
Clearly, both sets are $S^1$-invariant. Denote by $\crit(\Psi_{a, K})$ the set of critical points of $\Psi_{a, K}$.
Because $\Psi_{a,K}$ is $S^1$-invariant, $S^1\cdot u$ becomes a critical orbit if $u\in \crit(\Psi_{a, K})$.
Note that by the condition (F), the number of critical orbits of $\Psi_{a, K}$
is finite. Hence as usual we can make the following definition.

{\bf Definition 2.1.} {\it Suppose $u$ is a nonzero critical point of $\Psi_{a, K}$, and $\Nn$ is an $S^1$-invariant
open neighborhood of $S^1\cdot u$ such that $\crit(\Psi_{a,K})\cap (\Lm_{a,K}(u)\cap \Nn) = S^1\cdot u$.
Then the $S^1$-critical module of $S^1\cdot u$ is defined by
$$ C_{S^1,\; q}(\Psi_{a, K}, \;S^1\cdot u)
=H_{q}((\Lambda_{a, K}(u)\cap\Nn)_{S^1},\; ((\Lambda_{a,K}(u)\setminus S^1\cdot u)\cap\Nn)_{S^1}). $$
Similarly, we define the $S^1$-critical module $C_{S^1,\; q}(F_{a, K}, \;S^1\cdot x)$ of $S^1\cdot x$
for $F_{a, K}$.}

We fix $a$ and let $u_K\neq 0$ be a critical point of $\Psi_{a, K}$ with multiplicity $\mul(u_K)=m$,
that is, $u_K$ corresponds to a closed characteristic $(\tau, y)\subset\Sigma$ with $(\tau, y)$
being $m$-iteration of
some prime closed characteristic. Precisely, we have $u_K=-J\dot x+Kx$ with $x$
being a solution of (\ref{2.1}) and $x=\rho y(\frac{\tau t}{T})$ with
$\frac{\vf_a^\prime(\rho)}{\rho}=\frac{\tau}{aT}$.
Moreover, $(\tau, y)$ is a closed characteristic on $\Sigma$ with minimal period $\frac{\sigma}{m}$.
For any $p\in\N$ satisfying $p\sigma<aT$, we choose $K$
such that $pK\notin \frac{2\pi}{T}\Z$, then the $p$th iteration $u_{pK}^p$ of $u_K$ is given by $-J\dot x^p+pKx^p$,
where $x^p$ is the unique solution of (\ref{2.1}) corresponding to $(p\tau, y)$
and is a critical point of $F_{a, pK}$, that
is, $u_{pK}^p$ is the critical point of $\Psi_{a, pK}$ corresponding to $x^p$.

\medskip

{\bf Lemma 2.2.} (cf. Proposition 4.2 and Remark 4.4 of \cite{LLW14} ) {\it If $u_{pK}^p$ is non-degenerate,
i.e., $\nu_{pK}(u_{pK}^p)=1$, let
$\bb(x^p)=(-1)^{i_{pK}(u_{pK}^p)-i_{K}(u_{K})}=(-1)^{i^v(x^p)-i^v(x)}$, then}
\bea C_{S^1,q-d(pK)+d(K)}(F_{a,K},S^1\cdot x^p)
&=& C_{S^1,q}(F_{a,pK},S^1\cdot x^p)=C_{S^1,q}(\Psi_{a,pK},S^1\cdot u^p_{pK}) \nn\\
&=& \left\{\matrix{
     \Q, &\quad {\it if}\;\; q=i_{pK}(u_{pK}^p),\;\;{\it and}\;\;\bb(x^p)=1, \cr
     0, &\quad {\it otherwise}. \cr}\right.  \lb{2.10}\eea

\medskip

{\bf Theorem 2.3.} (cf. Theorem 1.1 of \cite{LLW14} and Theorem 1.2 of \cite{Vit2}) {\it Suppose that
$\Sg\in\H_{st}(2n)$ satisfying $^\#\T(\Sg)<+\infty$. Denote by $\{(\tau_j,y_j)\}_{1\le j\le k}$ all the
geometrically distinct prime closed characteristics. Then the following identities hold
\bea \sum_{1\le j\le k \atop \hat{i}(y_j)>0}\frac{\hat{\chi}(y_j)}{\hat{i}(y_j)}=\frac{1}{2},\qquad
\sum_{1\le j\le k \atop \hat{i}(y_j)<0}\frac{\hat{\chi}(y_j)}{\hat{i}(y_j)}=0,\lb{2.11}\eea
where $\hat{\chi}(y)\in\Q$ is the average Euler characteristic given by Definition 4.8 and Remark 4.9 of \cite{LLW14}.

In particular, if all $y^m$'s are non-degenerate for $m\ge 1$, then
\bea \hat{\chi}(y)=\left\{\matrix{
     (-1)^{i(y)}, &\quad {\it if}\;\; i(y^2)-i(y)\in 2\Z, \cr
     \frac{(-1)^{i(y)}}{2}, &\quad {\it otherwise}. \cr}\right.  \lb{2.12}\eea}

Let $F_{a, K}$ be a functional defined by (\ref{2.4}) for some $a, K\in\R$ large enough and let $\ep>0$ be
small enough such that $[-\ep, 0)$ contains no critical values of $F_{a, K}$. For $b$ large enough,
The normalized Morse series of $F_{a, K}$ in $ X^{-\ep}\setminus X^{-b}$
is defined, as usual, by
\be  M_a(t)=\sum_{q\ge 0,\;1\le j\le p} \dim C_{S^1,\;q}(F_{a, K}, \;S^1\cdot v_j)t^{q-d(K)},  \lb{2.13}\ee
where we denote by $\{S^1\cdot v_1, \ldots, S^1\cdot v_p\}$ the critical orbits of $F_{a, K}$ with critical
values less than $-\ep$. The Poincar\'e series of $H_{S^1, *}( X, X^{-\ep})$ is $t^{d(K)}Q_a(t)$, according
to Theorem 5.1 of \cite{LLW14}, if we set $Q_a(t)=\sum_{k\in \Z}{q_kt^k}$, then
$$   q_k=0 \qquad\qquad \forall\;k\in \mathring {I},  $$
where $I$ is an interval of $\Z$ such that $I \cap [i(\tau, y), i(\tau, y)+\nu(\tau, y)-1]=\emptyset$ for all
closed characteristics $(\tau,\, y)$ on $\Sigma$ with $\tau\ge aT$. Then by Section 6 of \cite{LLW14}, we have
$$  M_a(t)-\frac{1}{1-t^2}+Q_a(t) = (1+t)U_a(t),   $$
where $U_a(t)=\sum_{i\in \Z}{u_it^i}$ is a Laurent series with nonnegative coefficients.
If there is no closed characteristic with $\hat{i}=0$, then
\be   M(t)-\frac{1}{1-t^2}=(1+t)U(t),    \lb{2.14}\ee
where $M(t)=\sum_{p\in \Z}{M_pt^p}$ denotes $M_a(t)$ as $a$ tends to infinity. In addition, we also denote by $b_p$ the coefficient of $t^p$ of $\frac{1}{1-t^2}=\sum_{p\in \Z}{b_pt^p}$, i.e. there holds $b_p=1$, $\forall\ p\in2\N_0$ and $b_p=0$, $\forall\ p\not\in2\N_0$.

\setcounter{figure}{0}
\setcounter{equation}{0}
\section{The enhanced common index jump theorem of symplectic paths}

In \cite{Lon2} of 1999, Y. Long established the basic normal form
decomposition of symplectic matrices. Based on this result he
further established the precise iteration formulae of indices of
symplectic paths in \cite{Lon3} of 2000.

As in \cite{Lon3}, denote by
\bea
N_1(\lm, b) &=& \left(\matrix{\lm & b\cr
                                0 & \lm\cr}\right), \qquad {\rm for\;}\lm=\pm 1, \; b\in\R, \lb{3.1}\\
D(\lm) &=& \left(\matrix{\lm & 0\cr
                      0 & \lm^{-1}\cr}\right), \qquad {\rm for\;}\lm\in\R\bs\{0, \pm 1\}, \lb{3.2}\\
R(\th) &=& \left(\matrix{\cos\th & -\sin\th \cr
                           \sin\th & \cos\th\cr}\right), \qquad {\rm for\;}\th\in (0,\pi)\cup (\pi,2\pi), \lb{3.3}\\
N_2(e^{\th\sqrt{-1}}, B) &=& \left(\matrix{ R(\th) & B \cr
                  0 & R(\th)\cr}\right), \qquad {\rm for\;}\th\in (0,\pi)\cup (\pi,2\pi)\;\; {\rm and}\; \nn\\
        && \qquad B=\left(\matrix{b_1 & b_2\cr
                                  b_3 & b_4\cr}\right)\; {\rm with}\; b_j\in\R, \;\;
                                         {\rm and}\;\; b_2\not= b_3. \lb{3.4}\eea
Here $N_2(e^{\th\sqrt{-1}}, B)$ is non-trivial if $(b_2-b_3)\sin\theta<0$, and trivial
if $(b_2-b_3)\sin\theta>0$.

As in \cite{Lon3}, the $\diamond$-sum (direct sum) of any two real matrices is defined by
$$ \left(\matrix{A_1 & B_1\cr C_1 & D_1\cr}\right)_{2i\times 2i}\diamond
      \left(\matrix{A_2 & B_2\cr C_2 & D_2\cr}\right)_{2j\times 2j}
=\left(\matrix{A_1 & 0 & B_1 & 0 \cr
                                   0 & A_2 & 0& B_2\cr
                                   C_1 & 0 & D_1 & 0 \cr
                                   0 & C_2 & 0 & D_2}\right). $$

For every $M\in\Sp(2n)$, the homotopy set $\Omega(M)$ of $M$ in $\Sp(2n)$ is defined by
$$ \Om(M)=\{N\in\Sp(2n)\,|\,\sg(N)\cap\U=\sg(M)\cap\U\equiv\Gamma\;\mbox{and}
                    \;\nu_{\om}(N)=\nu_{\om}(M)\, \forall\om\in\Gamma\}, $$
where $\sg(M)$ denotes the spectrum of $M$,
$\nu_{\om}(M)\equiv\dim_{\C}\ker_{\C}(M-\om I)$ for $\om\in\U$.
The component $\Om^0(M)$ of $P$ in $\Sp(2n)$ is defined by
the path connected component of $\Om(M)$ containing $M$.

\medskip

{\bf Lemma 3.1.} (cf. \cite{Lon3}, Lemma 9.1.5 and List 9.1.12 of \cite{Lon4})
{\it For $M\in\Sp(2n)$ and $\om\in\U$, the splitting number $S_M^\pm(\om)$ (cf. Definition 9.1.4 of \cite{Lon4}) satisfies
\begin{eqnarray}
S_M^{\pm}(\om) &=& 0, \qquad {\it if}\;\;\om\not\in\sg(M).  \lb{3.5}\\
S_{N_1(1,a)}^+(1) &=& \left\{\matrix{1, &\quad {\rm if}\;\; a\ge 0, \cr
0, &\quad {\rm if}\;\; a< 0. \cr}\right. \lb{3.6}\eea

For any $M_i\in\Sp(2n_i)$ with $i=0$ and $1$, there holds }
\be S^{\pm}_{M_0\diamond M_1}(\om) = S^{\pm}_{M_0}(\om) + S^{\pm}_{M_1}(\om),
    \qquad \forall\;\om\in\U. \lb{3.7}\ee

We have the following decomposition theorem

\medskip

{\bf Theorem 3.2.} (cf. \cite{Lon3} and Theorem 1.8.10 of \cite{Lon4}) {\it For
any $M\in\Sp(2n)$, there is a path $f:[0,1]\to\Om^0(M)$ such that $f(0)=M$ and
\be f(1) = M_1\diamond\cdots\diamond M_k,  \lb{3.8}\ee
where each $M_i$ is a basic normal form listed in (\ref{3.1})-(\ref{3.4})
for $1\leq i\leq k$.}

\medskip

For every $\ga\in\mathcal{P}_\tau(2n)\equiv\{\ga\in C([0,\tau],Sp(2n))\ |\ \ga(0)=I_{2n}\}$, we extend
$\ga(t)$ to $t\in [0,m\tau]$ for every $m\in\N$ by
\bea \ga^m(t)=\ga(t-j\tau)\ga(\tau)^j \qquad \forall\;j\tau\le t\le (j+1)\tau \;\;
               {\rm and}\;\;j=0, 1, \ldots, m-1, \lb{3.9}\eea
as in P.114 of \cite{Lon2}. As in \cite{LoZ} and \cite{Lon4}, we denote the Maslov-type indices of
$\ga^m$ by $(i(\ga,m),\nu(\ga,m))$.

Then the following iteration formula from \cite{LoZ} and \cite{Lon4} can be obtained.

\medskip

{\bf Theorem 3.3.} (cf. Theorem 9.3.1 of \cite{Lon4}) {\it For any path $\ga\in\mathcal{P}_\tau(2n)$,
let $M=\ga(\tau)$ and $C(M)=\sum_{0<\th<2\pi}S_M^-(e^{\sqrt{-1}\th})$. We extend $\ga$ to $[0,+\infty)$
by its iterates. Then for any $m\in\N$ we have
\bea i(\ga,m)
&=& m(i(\ga,1)+S^+_{M}(1)-C(M))\nn\\
& & + 2\sum_{\th\in(0,2\pi)}E\left(\frac{m\th}{2\pi}\right)S^-_{M}(e^{\sqrt{-1}\th}) - (S_M^+(1)+C(M)), \lb{3.10}\eea
and
\be \hat{i}(\ga,1) = i(\ga,1) + S^+_{M}(1) - C(M) + \sum_{\th\in(0,2\pi)}\frac{\th}{\pi}S^-_{M}(e^{\sqrt{-1}\th}). \lb{3.11}\ee}

\medskip

The common index jump theorem (cf. Theorem 4.3 of \cite{LoZ}) for symplectic paths established by Long
and Zhu in 2002 has become one of the main tools to study the multiplicity and stability problems of
closed solution orbits in Hamiltonian and symplectic dynamics. Recently, the following enhanced version
of it has been obtained by Duan, Long and Wang in \cite{DLW}, which will play an important role in the
proofs in Section 4.

\medskip

{\bf Theorem 3.4.} (cf. Theorem 3.5 of \cite{DLW}) ({\bf The enhanced common index jump theorem for
symplectic paths}) {\it Let $\gamma_k\in\mathcal{P}_{\tau_k}(2n)$ for $k=1,\cdots,q$ be a finite
collection of symplectic paths. Let $M_k=\ga(\tau_k)$. We extend $\ga_k$ to $[0,+\infty)$ by (\ref{3.9})
inductively. Suppose
\be  \hat{i}(\ga_k,1) > 0, \qquad \forall\ k=1,\cdots,q.  \lb{3.12}\ee
Then for every integer $\bar{m}\in \N$, there exist infinitely many $(q+1)$-tuples
$(N, m_1,\cdots,m_q) \in \N^{q+1}$ such that for all $1\le k\le q$ and $1\le m\le \bar{m}$, there holds
\bea
\nu(\ga_k,2m_k-m) &=& \nu(\ga_k,2m_k+m) = \nu(\ga_k, m),   \lb{3.13}\\
i(\ga_k,2m_k+m) &=& 2N+i(\ga_k,m),                         \lb{3.14}\\
i(\ga_k,2m_k-m) &=& 2N-i(\ga_k,m)-2(S^+_{M_k}(1)+Q_k(m)),  \lb{3.15}\\
i(\ga_k, 2m_k)&=& 2N -(S^+_{M_k}(1)+C(M_k)-2\Delta_k),     \lb{3.16}\eea
where
\be \Delta_k = \sum_{0<\{m_k\th/\pi\}<\delta}S^-_{M_k}(e^{\sqrt{-1}\th}),\qquad
 Q_k(m) = \sum_{e^{\sqrt{-1}\th}\in\sg(M_k),\atop \{\frac{m_k\th}{\pi}\}
                   = \{\frac{m\th}{2\pi}\}=0}S^-_{M_k}(e^{\sqrt{-1}\th}). \lb{3.17}\ee
More precisely, by (4.10), (4.40) and (4.41) in \cite{LoZ} , we have
\bea m_k=\left(\left[\frac{N}{M\hat i(\gamma_k, 1)}\right]+\chi_k\right)M,\quad 1\le k\le q,\lb{3.18}\eea
where $\chi_k=0$ or $1$ for $1\le k\le q$ and $\frac{M\theta}{\pi}\in\Z$
whenever $e^{\sqrt{-1}\theta}\in\sigma(M_k)$ and $\frac{\theta}{\pi}\in\Q$
for some $1\le k\le q$.  Furthermore, by (4.20) in Theorem 4.1 of \cite{LoZ},
for any $\epsilon>0$, we can choose $N$ and $\{\chi_k\}_{1\le k\le q}$ such that}
\bea \left|\left\{\frac{N}{M\hat i(\gamma_k, 1)}\right\}-\chi_k\right|<\epsilon,\quad 1\le k\le q.\lb{3.19}\eea

\medskip

{\bf Theorem 3.5.} (cf. Theorem 2.1 of \cite{HuL} and Theorem 6.1 of \cite{LLo2}) {\it Suppose $\Sg\in \H_{st}(2n)$ and
$(\tau,y)\in \T(\Sigma)$. Then we have
\be i(y^m)\equiv i(m\tau,y)=i(y, m)-n,\quad \nu(y^m)\equiv\nu(m\tau, y)=\nu(y, m),
       \qquad \forall m\in\N, \lb{3.20}\ee
where $i(y^m)$ and $\nu(y^m)$ are the index and nullity of $(m\tau,y)$ defined in Section 2, $i(y, m)$ and $\nu(y, m)$
are the Maslov-type index and nullity of $(m\tau,y)$ (cf. Section 5.4 of \cite{Lon3}). In particular, we have
$\hat{i}(\tau,y)=\hat{i}(y,1)$, where $\hat{i}(\tau ,y)$ is given in Section 2, $\hat{i}(y,1)$
is the mean Maslov-type index (cf. Definition 8.1 of \cite{Lon4}). Hence we denote it simply by $\hat{i}(y)$.}

\setcounter{figure}{0}
\setcounter{equation}{0}
\section{Proof of Theorem 1.2}

In order to prove Theorem 1.2, we make the following assumptions

\medskip

{\bf (FCHe)} {\it Suppose that on  a non-degenerate compact star-shaped hypersurface $\Sigma$
in $\R^{2n}$ with $n$ {\bf even}, there exist only finitely many prime closed characteristics
$\{(\tau_k,y_k)\}_{k=1}^q$ with $\hat{i}(y_k,1)>0$ and $i(y_k,m)\neq -1$ for all $m\in\N$.}

\medskip

{\bf (FCHo)} {\it Suppose that on  a non-degenerate compact star-shaped hypersurface $\Sigma$
in $\R^{2n}$ with $n$ {\bf odd}, there exist only finitely many prime closed characteristics
$\{(\tau_k,y_k)\}_{k=1}^q$ with $\hat{i}(y_k,1)>0$ and $i(y_k,m)\notin\{-2,-1,0\}$ for all $m\in\N$.}

\medskip

Since $\Sigma$ is index perfect, then the assumptions (FCHe) and (FCHo) hold by Definition 1.1.
We denote by $\ga_k\equiv \ga_{y_k}$ the associated symplectic path of $(\tau_k,y_k)$ for $1\le k\le q$.
Then by Lemma 3.3 of \cite{HuL} and Lemma 3.2 of \cite{Lon1},
there exists $P_k\in Sp(2n)$ and $U_k\in Sp(2n-2)$ such that
\bea M_k\equiv\ga_k(\tau_k)=P_k^{-1}(N_1(1,1)\dm U_k)P_k,\qquad 1\le k\le q,\lb{4.1}\eea
where every $U_k$ has the following form:
\bea
&& R(\th_1)\,\dm\,\cdots\,\dm\,R(\th_r)\,\dm\,D(\pm 2)^{\dm s} \nn\\
&& \dm\,N_2(e^{\aa_{1}\sqrt{-1}},A_{1})\,\dm\,\cdots\,\dm\,N_2(e^{\aa_{r_{\ast}}\sqrt{-1}},A_{r_{\ast}})
  \,\dm\,N_2(e^{\bb_{1}\sqrt{-1}},B_{1})\,\dm\,\cdots\,\dm\,N_2(e^{\bb_{r_{0}}\sqrt{-1}},B_{r_{0}}), \nn\eea
where $\frac{\th_{j}}{2\pi}\not\in\Q$ for $1\le j\le r$; $\frac{\aa_{j}}{2\pi}\not\in\Q$ for $1\le j\le r_{\ast}$;
$\frac{\bb_{j}}{2\pi}\not\in\Q$ for $1\le j\le r_0$ and
\be r+ s +2r_{\ast} + 2r_0 = n-1. \lb{4.2}\ee

\medskip

{\bf Proof of Theorem 1.2.}

\medskip

We prove Theorem 1.2 in two cases:

\medskip

{\bf Case 1.} {\it $n$ is even.}

\medskip

By the assumption (FCHe), we have $\hat{i}(y_k)=\hat{i}(y_k,1)>0, 1\le k\le q$,
which implies that $i(y_k,m)\rightarrow +\infty$ as $m\rightarrow +\infty$.  So the positive integer
$\bar{m}$ defined by
\bea \bar{m}=\max_{1\le k\le q}\left\{\min\{m_0\in\N\ |\ i(y_k,m+l)\ge i(y_k,l)+n+1, \quad
    \forall\ l\ge 1, m\geq m_0\}\right\} \lb{4.3}\eea
is well-defined and finite.

For the integer $\bar{m}$ defined in (\ref{4.3}), it follows from Theorem 3.4
that there exist infinitely many $q+1$-tuples $(N, m_1, \cdots, m_q)\in\N^{q+1}$ such that for any
$1\le k\le q$, there holds
\bea
\bar{m}+2 &\le& \min\{2m_k,\ 1\le k\le q\},\lb{4.4}\\
i(y_k,{2m_k-m}) &=& 2N-2-i(y_k,m),\quad 1\le m\le\bar{m}, \lb{4.5}\\
i(y_k,{2m_k}) &=& 2N-1-C(M_k)+2\Delta_k,\lb{4.6}\\
i(y_k,{2m_k+m})&=& 2N+i(y_k,m),\quad 1\le m\le\bar{m},\lb{4.7}\eea
where, note that $S^+_{M_k}(1)=1,Q_k(m)=0$, $\forall\ m\ge 1$ by (\ref{4.1})-(\ref{4.2}).

By the definition (\ref{4.3}) of $\bar{m}$ and (\ref{4.6}), for $m\ge \bar{m}+1$, we obtain
\bea
&& i(y_k,2m_k-m)\le i(y_k,2m_k)-n-1=2N-n-2+2\Delta_k-C(M_k) \le 2N-3,\lb{4.8}\\
&& i(y_k,2m_k+m)\ge i(y_k,2m_k)+n+1=2N+n-C(M_k)+2\Delta_k \ge 2N+1,\lb{4.9}\eea
where we use the facts $2\Delta_k-C(M_k)\le n-1$ and $C(M_k)\le n-1$.

Then by (\ref{4.5})-(\ref{4.9}) and Theorem 3.5, we obtain
\bea  i(y_k^{2m_k-m})&\le& 2N-n-3,\qquad \forall\ \bar{m}+1\le m\le 2m_k-1,\lb{4.10} \\
 i(y_k^{2m_k-m}) &=& 2N-2n-2-i(y_k^m),\qquad \forall\ 1\le m\le \bar{m}, \lb{4.11}\\
 i(y_k^{2m_k}) &=& 2N-C(M_k)+2\Delta_k-n-1,\lb{4.12}\\
 i(y_k^{2m_k+m}) &=& 2N+i(y_k^m),\qquad \forall\ 1\le m\le \bar{m},\lb{4.13}\\
 i(y_k^{2m_k+m}) &\ge& 2N-n+1,\qquad \forall\ m\ge \bar{m}+1.\lb{4.14} \eea

\medskip

{\bf Claim 1:} {\it For $N\in\N$ in Theorem 3.4 satisfying (\ref{4.10})-(\ref{4.14}), we have
\bea \sum_{1\le k\le q} 2m_k\hat{\chi}(y_k)=N.\lb{4.15}\eea}

In fact, let $\ep<\frac{1}{1+2M\sum_{1\le k\le q}|\hat{\chi}(y_k)|}$, by Theorem 2.3 and
(\ref{3.18})-(\ref{3.19}), it yields
\bea \left|N-\sum_{k=1}^q 2m_k\hat{\chi}(y_k)\right|
&=& \left|\sum_{k=1}^q\frac{2N\hat{\chi}(y_k)}{\hat{i}(y_k)}-\sum_{k=1}^q 2\hat{\chi}(y_k)
    \left(\left[\frac{N}{M\hat{i}(y_k)}\right]+\chi_k\right)M\right|\nn\\
&\le& 2M\sum_{k=1}^q |\hat{\chi}(y_k)|\left|\left\{\frac{N}{M\hat{i}(y_k)}\right\}-\chi_k\right|.\nn\\
&<& 2M\ep\sum_{k=1}^q|\hat{\chi}(y_k)|\nn\\
&<& 1,\lb{4.16}\eea
which proves Claim 1 since each $2m_k\hat{\chi}(y_k)$ is an integer.

\medskip

Now by Lemma 2.2, it yields
\bea
&& \sum_{m=1}^{2m_k} (-1)^{d(K)+i(y_k^m)}\dim C_{S^1,d(K)+i(y_k^{m})}(F_{a, K},S^1\cdot x_k^m)\nn\\
&& = \sum_{m=1}^{2m_k} (-1)^{i(y_k^m)}\dim C_{S^1,d(K)+i(y_k^{m})}(F_{a, K},S^1\cdot x_k^m)\nn\\
&& = \sum_{i=0}^{m_k-1} \sum_{m=2i+1}^{2i+2} (-1)^{i(y_k^m)} \dim C_{S^1,d(K)+i(y_k^{m})}(F_{a, K},S^1\cdot x_k^m)\nn\\
&& = \sum_{i=0}^{m_k-1} \sum_{m=1}^{2} (-1)^{i(y_k^m)} \dim C_{S^1,d(K)+i(y_k^{m})}(F_{a, K},S^1\cdot x_k^m)\nn\\
&& = m_k \sum_{m=1}^{2} (-1)^{i(y_k^m)} \dim C_{S^1,d(K)+i(y_k^{m})}(F_{a, K},S^1\cdot x_k^m)\nn\\
&& = 2m_k\hat{\chi}(y_k),\qquad \forall\ 1\le k\le q, \lb{4.17}\eea
where $x_k$ is the critical point of $F_{a, K}$ corresponding to $y_k$, $d(K)=2n([KT/{2\pi}]+1)$.

Note that the set $\{i(y_k^m): m\in\N, 1\le k\le q\}$ is bounded below since the mean index $\hat{i}(y_k)>0$ for any $1\le k\le q$.
So the alternative sum $\sum_{p=-\infty}^{2N-n-1}(-1)^p M_p$ is a finite sum by Lemma 2.2, where the Morse-type number $M_p$ is defined in (\ref{2.14}).

For $1\le k\le q$, by (\ref{4.14}) and Lemma 2.2, we know that all $y_k^{2m_k+m}$'s with $m\ge\bar{m}+1$
have no contribution to the alternative sum $\sum_{p=-\infty}^{2N-n-1}(-1)^p M_p$. Similarly again by
Lemma 2.2 and (\ref{4.10}), all $y_k^{2m_k-m}$'s with $\bar{m}+1\le m\le 2m_k-1$ only have contribution to
$\sum_{p=-\infty}^{2N-n-1}(-1)^p M_p$.

For $1\le k\le q$, by (\ref{4.13}) and Lemma 2.2, we know that all $y_k^{2m_k+m}$'s with $-n\le i(y_k^m)$
and $1\le m\le\bar{m}$ have no contribution to the alternative sum $\sum_{p=-\infty}^{2N-n-1}(-1)^p M_p$.
Similarly again by Lemma 2.2 and (\ref{4.11}), all $y_k^{2m_k-m}$'s with $-n\le i(y_k^m)$ and
$1\le m\le\bar{m}$ only have contribution to $\sum_{p=-\infty}^{2N-n-1}(-1)^p M_p$.

Since $i(y_k^{m})\neq -n-1$ for any $m\ge 1$ by the assumption (FCHe) and Theorem 3.5, for $1\le k\le q$,
we set
\bea
M_+^e(k) &=& ^{\#}\{1\le m\le \bar{m}\ |\ i(y_k^{m})\le -n-2,\ i(y_k^{2m_k+m})-i(y_k)\in 2\Z,\ i(y_k)\in 2\Z\},\nn\\
M_+^o(k) &=& ^{\#}\{1\le m\le \bar{m}\ |\ i(y_k^{m})\le -n-2,\ i(y_k^{2m_k+m})-i(y_k)\in 2\Z,\ i(y_k)\in 2\Z-1\},\nn\\
M_-^e(k) &=& ^{\#}\{1\le m\le \bar{m}\ |\ i(y_k^{m})\le -n-2,\ i(y_k^{2m_k-m})-i(y_k)\in 2\Z,\ i(y_k)\in 2\Z\},\nn\\
M_-^o(k) &=& ^{\#}\{1\le m\le \bar{m}\ |\ i(y_k^{m})\le -n-2,\ i(y_k^{2m_k-m})-i(y_k)\in 2\Z,\ i(y_k)\in 2\Z-1\},\nn\eea
which, together with $i(y_k^{2m_k+m})-i(y_k^{2m_k-m})\in 2\Z$ by (\ref{4.11}) and (\ref{4.13}), yields
\be M_+^{e}(k)=M_-^{e}(k),\quad M_+^{o}(k)=M_-^{o}(k),\qquad 1\le k\le q.\lb{4.18}\ee

For $1\le k\le q$ and $1\le m\le \bar{m}$ satisfying $i(y_k^{m})\le -n-2$, by (\ref{4.11}) and (\ref{4.13})
it yields
\be i(y_k^{2m_k-m})\ge 2N-n,\qquad i(y_k^{2m_k+m})\le 2N-n-2.\lb{4.19}\ee

So, for $1\le k\le q$, by (\ref{4.19}) and Lemma 2.2, we know that all $y_k^{2m_k+m}$'s with $i(y_k^m)\le -n-2$
and $1\le m\le\bar{m}$ only have contribution to the alternative sum $\sum_{p=-\infty}^{2N-n-1}(-1)^p M_p$, and
all $y_k^{2m_k-m}$'s with $i(y_k^m)\le -n-2$ and $1\le m\le\bar{m}$ have no contribution to
$\sum_{p=-\infty}^{2N-n-1}(-1)^p M_p$.

Thus for the Morse-type numbers $M_p$'s in (\ref{2.14}), by (\ref{4.17})-(\ref{4.19}) we have
\bea \sum_{p=-\infty}^{2N-n-1}(-1)^p M_p
&=& \sum_{k=1}^{q}\ \sum_{1\le m\le 2m_k+\bar{m} \atop i(y_k^{m})\le 2N-n-1}
    (-1)^{d(K)+i(y_k^m)}\dim C_{S^1,d(K)+i(y_k^{m})}(F_{a, K},S^1\cdot x_k^m)\nn\\
&=& \sum_{k=1}^{q}\ \sum_{m=1}^{2m_k} (-1)^{d(K)+i(y_k^m)}\dim C_{S^1,d(K)+i(y_k^{m})}(F_{a, K},S^1\cdot x_k^m)\nn\\
& & +\sum_{k=1}^q \left[M_+^e(k)-M_+^o(k)\right]-\sum_{k=1}^q \left[M_-^e(k)-M_-^o(k)\right]\nn\\
& & -\sum_{1\le k\le q \atop i(y_k^{2m_k})\ge 2N-n} (-1)^{i(y_k^{2m_k})}
              \dim C_{S^1,d(K)+i(y_k^{2m_k})}(F_{a, K},S^1\cdot x_k^{2m_k})\nn\\
&=& \sum_{k=1}^{q} 2m_k\hat{\chi}(y_k)\nn\\
& & -\sum_{1\le k\le q \atop i(y_k^{2m_k})\ge 2N-n} (-1)^{i(y_k^{2m_k})}
             \dim C_{S^1,d(K)+i(y_k^{2m_k})}(F_{a, K},S^1\cdot x_k^{2m_k}). \lb{4.20}\eea

In order to exactly know whether the iterate $y_k^{2m_k}$ of $y_k$ has contribution to the alternative sum
$\sum_{p=-\infty}^{2N-n-1}(-1)^p M_p(k)$,  $1\le k\le q$, we set
\bea
N_+^e &=& ^{\#}\{1\le k\le q\ |\ i(y_k^{2m_k})\ge 2N-n,\ i(y_k^{2m_k})-i(y_k)\in 2\Z,\ i(y_k)\in 2\Z\},\lb{4.21}\\
N_+^o &=& ^{\#}\{1\le k\le q\ |\ i(y_k^{2m_k})\ge 2N-n,\ i(y_k^{2m_k})-i(y_k)\in 2\Z,\ i(y_k)\in 2\Z-1\},\lb{4.22}\\
N_-^e &=& ^{\#}\{1\le k\le q\ |\ i(y_k^{2m_k})\le 2N-n-2,\ i(y_k^{2m_k})-i(y_k)\in 2\Z,\ i(y_k)\in 2\Z\},\lb{4.23}\\
N_-^o &=& ^{\#}\{1\le k\le q\ |\ i(y_k^{2m_k})\le 2N-n-2,\ i(y_k^{2m_k})-i(y_k)\in 2\Z,\ i(y_k)\in 2\Z-1\}.\lb{4.24}\eea

Note that $\sum_{p=-\infty}^{2N-n-1}(-1)^p M_p$ is a finite sum by the positive mean index assumption. Thus by Claim 1, (\ref{4.20}), the definitions of $N^{e}_+$ and $N_+^o$ and (\ref{2.14}), we have
\bea N+N_+^o-N_+^e
&=& \sum_{k=1}^q 2m_k\hat{\chi}(y_k)+N_+^o-N_+^e  \nn\\
&=& \sum_{p=-\infty}^{2N-n-1}(-1)^p M_p \nn\\
&\le& \sum_{p=-\infty}^{2N-n-1}(-1)^p b_p=\sum_{p=0}^{2N-n-2}b_p \nn\\
&=& N-\frac{n}{2},\lb{4.25}\eea
which implies
\be N_+^e\ge \frac{n}{2}.\lb{4.26}\ee

Similar to (\ref{4.10})-(\ref{4.14}), it follows from Theorem 3.4 that there exist also
infinitely many $(q+1)$-tuples $(N', m_1', \cdots, m_q')\in\N^{q+1}$ such that for any $1\le k\le q$,
there holds
\bea
i(y_k^{2m_k'-m})&\le& 2N'-n-3,\qquad \forall\ \bar{m}+1\le m\le 2m_k'-1,\lb{4.27} \\
i(y_k^{2m_k'-m}) &=& 2N'-2n-2-i(y_k^m),\qquad \forall\ 1\le m\le \bar{m}, \lb{4.28}\\
i(y_k^{2m_k'}) &=& 2N'-C(M_k)+2\Delta_k'-n-1,\lb{4.29}\\
i(y_k^{2m_k'+m}) &=& 2N'+i(y_k^m),\qquad \forall\ 1\le m\le \bar{m},\lb{4.30}\\
i(y_k^{2m_k'+m}) &\ge& 2N'-n+1,\qquad \forall\ m\ge \bar{m}+1.\lb{4.31} \eea
where, furthermore, $\Delta_k$ and $\Delta_k'$ satisfy the following relationship
\bea \Delta_k' + \Delta_k = C(M_k),\qquad\forall\ 1\le k\le q,\lb{4.32}\eea
which exactly follows from (especially the term (c) of) Theorem 4.2 of \cite{LoZ}.

Similarly, we define
\bea
N_+^{'e} &=& ^{\#}\{1\le k\le q\ |\ i(y_k^{2m_k'})\ge 2N'-n,\ i(y_k^{2m_k'})-i(y_k)\in 2\Z_0,\ i(y_k)\in 2\Z\},   \lb{4.33}\\
N_+^{'o} &=& ^{\#}\{1\le k\le q\ |\ i(y_k^{2m_k'})\ge 2N'-n,\ i(y_k^{2m_k'})-i(y_k)\in 2\Z_0,\ i(y_k)\in 2\Z-1\}, \lb{4.34}\\
N_-^{'e} &=& ^{\#}\{1\le k\le q\ |\ i(y_k^{2m_k'})\le 2N'-n-2,\ i(y_k^{2m_k'})-i(y_k)\in 2\Z_0,\ i(y_k)\in 2\Z\},   \lb{4.35}\\
N_-^{'o} &=& ^{\#}\{1\le k\le q\ |\ i(y_k^{2m_k'})\le 2N'-n-2,\ i(y_k^{2m_k'})-i(y_k)\in 2\Z_0,\ i(y_k)\in 2\Z-1\}. \lb{4.36}\eea

So by (\ref{4.29}) and (\ref{4.32}) it yields
\be i(y_k^{2m_k'}) = 2N'-C(M_k)+2(C(M_k)-\Delta_k)-n-1=2N'+C(M_k)-2\Delta_k-n-1. \lb{4.37}\ee

So by definitions (\ref{4.21})-(\ref{4.24}) and (\ref{4.33})-(\ref{4.36}) we have
\be N_{\pm}^e=N_{\mp}^{'e},\qquad N_{\pm}^o=N_{\mp}^{'o}. \lb{4.38}\ee

Thus, through carrying out the similar arguments to (\ref{4.25})-(\ref{4.26}), by Claim 1, the definitions of
$N^{'e}_+$ and $N^{'o}_+$  and (\ref{2.14}), we have
\bea N'+N_+^{'o}-N_+^{'e}&=&\sum_{k=1}^q 2m_k\hat{\chi}(y_k)+N_+^{'o}-N_+^{'e}\nn\\
&=&\sum_{p=-\infty}^{2N'-n-1}(-1)^p M_p\nn\\
&\le&\sum_{p=-\infty}^{2N'-n-1}(-1)^p b_p=\sum_{p=0}^{2N'-n-2}b_p\nn\\
&=& N'-\frac{n}{2}\lb{4.39}\eea
which, together with (\ref{4.38}), implies
\bea N_-^{e}=N_+^{'e}\ge\frac{n}{2}.\lb{4.40}\eea

So by (\ref{4.26}) and (\ref{4.40}) it yields
\bea q\ge N_+^{e}+N_-^{e}\ge n. \lb{4.41}\eea

In addition, any hyperbolic closed characteristic $y_k$ must have $i(y_k^{2m_k})=2N-n-1$ since there holds
$C(M_k)=0$ in the hyperbolic case. However, by (\ref{4.21}) and (\ref{4.23}), there exist at least
$(N_+^{o}+N_-^{o})$ closed characteristic with even indices $i(y_k^{2m_k})$. So all these
$(N_+^{o}+N_-^{o})$ closed characteristics are non-hyperbolic. Then (\ref{4.41}) shows that there exist
at least $n$ distinct non-hyperbolic closed characteristics. And (\ref{4.21}), (\ref{4.23}) and (\ref{4.41})
show that all these non-hyperbolic closed characteristics and their iterations have even Maslov-type indices.
This completes the proof of Case 1.

\medskip

{\bf Case 2.} {\it $n$ is odd.}

\medskip

In this case, the assumption (FCHo) holds. Here the arguments are similar to those in the proof of Case 1. So we
only give some different parts in the proof and omit some details.

\medskip

{\bf Claim 2:} {\it There exist at least $(n-1)$ geometrically distinct non-hyperbolic closed characteristics
denoted by $\{y_k\}_{k=1}^{n-1}$ with odd Maslov-type indices on such hypersurface $\Sg$.}

\medskip

Here one crucial and different point from the proof of Case 1 is that we need to consider the alternative sum
$\sum_{p=-\infty}^{2N-n}(-1)^p M_p$ (cf. (\ref{4.44})) instead of the one $\sum_{p=-\infty}^{2N-n-1}(-1)^p M_p$
(cf. (\ref{4.20})). This difference is mainly due to the different parity of $n$. Since the method is similar
to that in proof of Case 1, we only list some necessary parts.

At first, there holds $i(y_k^m)\notin\{ -n-2,-n-1,-n\}$ for any $m\ge 1$ by the assumption (FCHo) and Theorem
3.5, for $1\le k\le q$, we set
\bea
\bar{M}_+^e(k) &=& ^{\#}\{1\le m\le \bar{m}\ |\ i(y_k^{m})\le -n-3,\ i(y_k^{2m_k+m})-i(y_k)\in 2\Z,\ i(y_k)\in 2\Z\},\nn\\
\bar{M}_+^o(k) &=& ^{\#}\{1\le m\le \bar{m}\ |\ i(y_k^{m})\le -n-3,\ i(y_k^{2m_k+m})-i(y_k)\in 2\Z,\ i(y_k)\in 2\Z-1\},\nn\\
\bar{M}_-^e(k) &=& ^{\#}\{1\le m\le \bar{m}\ |\ i(y_k^{m})\le -n-3,\ i(y_k^{2m_k-m})-i(y_k)\in 2\Z,\ i(y_k)\in 2\Z\},\nn\\
\bar{M}_-^o(k) &=& ^{\#}\{1\le m\le \bar{m}\ |\ i(y_k^{m})\le -n-3,\ i(y_k^{2m_k-m})-i(y_k)\in 2\Z,\ i(y_k)\in 2\Z-1\},\nn\eea
which, together with $i(y_k^{2m_k+m})-i(y_k^{2m_k-m})\in 2\Z$ by (\ref{4.11}) and (\ref{4.13}), yields
\bea \bar{M}_+^{e}(k)=\bar{M}_-^{e}(k),\quad \bar{M}_+^{o}(k)=\bar{M}_-^{o}(k),\qquad 1\le k\le q.\lb{4.42}\eea

For $1\le k\le q$ and $1\le m\le \bar{m}$ satisfying $i(y_k^{m})\le -n-3$, by (\ref{4.11}) and (\ref{4.13}) it yields
\bea i(y_k^{2m_k-m})\ge 2N-n+1,\qquad i(y_k^{2m_k+m})\le 2N-n-3.\lb{4.43}\eea

Then, similarly to the equation (\ref{4.20}), we have
\bea \sum_{p=-\infty}^{2N-n}(-1)^p M_p
&=& \sum_{k=1}^{q} 2m_k\hat{\chi}(y_k)\nn\\
& & -\sum_{1\le k\le q \atop i(y_k^{2m_k}) \ge 2N-n+1} (-1)^{i(y_k^{2m_k})}
      \dim C_{S^1,d(K)+i(y_k^{2m_k})}(F_{a, K},S^1\cdot x_k^{2m_k}). \lb{4.44}\eea

Denote by $H_+^e,H_+^o,H_-^e,H_-^o$ the numbers similarly defined by (\ref{4.21})-(\ref{4.24}) where $2N-n$ and
$2N-n-2$ are replaced by $2N-n+1$ and $2N-n-3$, respectively.

Then by Claim 1, (\ref{4.44}), the definitions of $H^{e}_+$ and $H_+^o$, and (\ref{2.14}), we have
\bea N+H_+^o-H_+^e&=&\sum_{k=1}^q 2m_k\hat{\chi}(y_k)+H_+^o-H_+^e\nn\\
&=&\sum_{p=-\infty}^{2N-n}(-1)^p M_p \nn\\
&=&-(M_{2N-n}-M_{2N-n-1}+\cdots-M_1+M_0-M_{-1}+\cdots)\nn\\
&\le&-(b_{2N-n}-b_{2N-n-1}+\cdots-b_1+b_0)\nn\\
&=&\sum_{p=0}^{2N-n-1} b_p=\frac{2N-n-1}{2}+1\nn\\
&=& N-\frac{n-1}{2},\lb{4.45}\eea
which yields
\bea H_+^e\ge H_+^e-H_+^o\ge \frac{n-1}{2}.\lb{4.46}\eea

Similarly, denote by $H_+^{'e},H_+^{'o},H_-^{'e},H_-^{'o}$ the numbers similarly defined by
(\ref{4.33})-(\ref{4.36}) where $2N'-n$ and $2N'-n-2$ are replaced by $2N'-n+1$ and $2N'-n-3$, respectively,
and these numbers satisfy the following relationship
\be H_{\pm}^e=H_{\mp}^{'e},\qquad H_{\pm}^o=H_{\mp}^{'o}. \lb{4.47}\ee

Similarly to the inequality (\ref{4.39}), by the same arguments above and (\ref{4.47}) we can obtain
\bea H_-^e=H_+^{'e}\ge H_+^{'e}-H_+^{'o}\ge \frac{n-1}{2}.\lb{4.48}\eea

Therefore it follows from (\ref{4.46}) and (\ref{4.48}) that
\bea q\ge H_+^e+H_-^e\ge \frac{n-1}{2}+\frac{n-1}{2}=n-1.\lb{4.49}\eea

By the same arguments in the proof of Case 1, it follows from the definitions of $H_+^e$ and $H_-^e$ that
these $(n-1)$ distinct closed geodesics are non-hyperbolic, and the Viterbo indices of them and their iterations
are even, and thus the Maslov-type indices of them and their iterations are odd. This completes the proof of
Claim 2.

\medskip

{\bf Claim 3:} {\it There exist at least another geometrically distinct closed characteristic different from
those found in Claim 2 with odd Maslov-type indices on such hypersurface $\Sg$.}

\medskip

In fact, for those $(n-1)$ distinct closed characteristics $\{y_k\}_{k=1}^{n-1}$ found in Claim 2, there holds
$i(y_k^{2m_k})\neq 2N-n-1$ by  the definitions of $H_+^e$ and $H_-^e$, which, together with
(\ref{4.10})-(\ref{4.14}) and the assumption (FCHo), yields
\bea i(y_k^m)\neq 2N-n-1,\qquad m\ge 1,\quad k=1,\cdots,n-1.\lb{4.50}\eea
Then by Lemma 2.2 it yields
\bea \sum_{1\le k\le n-1 \atop m\ge 1} \dim C_{S^1,d(K)+2N-n-1}(F_{a, K},S^1\cdot x_k^m)=0.\lb{4.51} \eea

By (\ref{4.10})-(\ref{4.11}), (\ref{4.13})-(\ref{4.14}) and the assumption (FCHo), it yields $i(y_k^m)\neq 2N-n-1$
for any $m\neq 2m_k$ and $1\le k\le q$. Therefore, by (\ref{4.51}) and (\ref{2.14}) we obtain
\bea \sum_{n\le k\le q}\dim C_{S^1,d(K)+2N-n-1}(F_{a, K},S^1\cdot x_k^{2m_k})
&=& \sum_{n\le k\le q,\ m\ge 1}\dim C_{S^1,d(K)+2N-n-1}(F_{a, K},S^1\cdot x_k^m) \nn\\
&=& \sum_{1\le k\le q,\ m\ge 1}\dim C_{S^1,d(K)+2N-n-1}(F_{a, K},S^1\cdot x_k^m)\nn\\
&=& M_{2N-n-1}\ge b_{2N-n-1}=1.\lb{4.52}\eea

Now by (\ref{4.52}) and Lemma 2.2, it yields that there exist at least another closed characteristic $y_n$
with $i(y_n^{2m_n})=2N-n-1$ and $i(y_n^{2m_n})-i(y_n)\in 2\Z$. Thus $y_n$ and its iterations have odd
Maslov-type indices. This completes the proof of Claim 3.

\medskip

Now for Case 2, Theorem 1.2 follows from Claim 2 and Claim 3. The proof of Theorem 1.2 is complete. \hfill\hb

\bibliographystyle{abbrv}

\end{document}